\documentclass{amsart}

\usepackage{amsfonts, amssymb, amsmath, eucal, verbatim, amsthm, amscd, enumerate}
\usepackage{graphicx}
\usepackage{setspace}

\newtheorem{theorem}{Theorem}

\theoremstyle{definition}

\theoremstyle{remark}
\newtheorem*{remark}{Remark}

\newtheorem*{acknowledgement}{Acknowledgement}

\numberwithin{equation}{section}
\begin{document}

\author{Jonathan Bennett, Neal Bez, Susana Guti\'errez}
\address{Jonathan Bennett, Neal Bez, Susana Guti\'errez: School of Mathematics, The Watson Building, University of Birmingham, Edgbaston,
Birmingham, B15 2TT, England} \email{[j.bennett,n.bez,s.gutierrez]@bham.ac.uk}
\author{Sanghyuk Lee}
\address{Sanghyuk Lee: Department of Mathematical Sciences, Seoul National University, Seoul 151-747, Korea}
\email{shklee@snu.ac.kr}
\title{On the Strichartz estimates for the kinetic transport equation}
\maketitle \thispagestyle{empty}

\begin{abstract}
We show that the endpoint Strichartz estimate for the kinetic transport equation
is false in all dimensions. We also present a new approach to proving the non-endpoint cases
using multilinear analysis.
\end{abstract}

\section{Introduction}

The solution of the kinetic transport equation 
\begin{equation*}
\partial_t f(t,x,v) + v \cdot \nabla_x f(t,x,v) = 0, \qquad f(0,x,v) = f^0(x,v) 
\end{equation*}
for $(t,x,v) \in \mathbb{R} \times \mathbb{R}^d \times \mathbb{R}^d$, satisfies the Strichartz estimates\footnote{We write $X \lesssim Y$ and $Y \gtrsim X$ if $X \leq CY$ for some finite constant $C$ depending at most on the parameters $(q,p,r,a,d)$, and $X \sim Y$ if $X \lesssim Y$ and $X \gtrsim Y$.}
\begin{equation} \label{e:Strichartz}
\|f\|_{L^q_tL^p_xL^r_v} \lesssim \|f^0\|_{L^a_{x,v}},
\end{equation}
where
\begin{equation} \label{e:sufficient}
\frac{2}{q} = d\bigg(\frac{1}{r} - \frac{1}{p}\bigg), \qquad \frac{1}{a} = \frac{1}{2}\bigg(\frac{1}{r} + \frac{1}{p}\bigg), \qquad q > a, \qquad p \geq a.
\end{equation}
With $(q,p,r,a)$ satisfying \eqref{e:sufficient}, but with the further condition $q > 2 \geq a$, this was proved by Castella and Perthame \cite{CastellaPerthame}, and it was observed by Keel and Tao \cite{KeelTao} that this latter condition can be relaxed to $q > a$ and hence \eqref{e:sufficient} suffices. In \cite{KeelTao} it is tentatively conjectured that the Strichartz estimate \eqref{e:Strichartz} holds at the endpoint $q=a$, at least for $d > 1$. Using the invariance under the transformations
\begin{equation} \label{e:transformations}
f^0 \leftrightarrow (f^0)^\lambda, \qquad f \leftrightarrow f^\lambda, \qquad (q,p,r,a) \leftrightarrow \bigg(\frac{q}{\lambda}, \frac{p}{\lambda}, \frac{r}{\lambda}, \frac{a}{\lambda}\bigg) 
\end{equation}
this conjectured endpoint can be (and usually is) stated for initial data in $L^2_{x,v}$ as
\begin{equation} \label{e:endpointL2}
\|f\|_{L^2_tL^{\frac{2d}{d-1}}_xL^{\frac{2d}{d+1}}_v} \lesssim \|f^0\|_{L^2_{x,v}}.
\end{equation}
The main purpose of this paper is to disprove this conjecture.
\begin{theorem} \label{t:fail}
The endpoint Strichartz estimate \eqref{e:endpointL2} for the kinetic transport equation fails for all $d \geq 1$.
\end{theorem}
The case $d=1$ of Theorem \ref{t:fail} was proved in \cite{GuoPeng} and \cite{OvcharovKakeya} by different arguments (where the norm in $x$ on the left-hand side is $L^\infty_x$). 

We remark that for the free Schr\"odinger propagator, the endpoint Strichartz estimate fails in the case $d = 2$ (see \cite{MontSmith}) but is true for all $d > 2$, as shown in the landmark paper of Keel and Tao \cite{KeelTao}. Thus, Theorem \ref{t:fail} highlights a fundamental difference in the Strichartz estimates for these related equations.

In the next section, we prove Theorem \ref{t:fail}. In the final section, we provide a new proof of the Strichartz estimates \eqref{e:Strichartz} in all non-endpoint cases using multilinear analysis.

\section{Proof of Theorem \ref{t:fail}}
Using the invariance under the transformation \eqref{e:transformations} with $\lambda = \frac{2d}{d+1}$,  estimate \eqref{e:endpointL2} is equivalent to
\begin{equation} \label{e:endpointtrans}
\|f\|_{L^{\frac{d+1}{d}}_tL^{\frac{d+1}{d-1}}_xL^{1}_v} \lesssim \|f^0\|_{L^{\frac{d+1}{d}}_{x,v}}.
\end{equation}
Since $f(t,x,v) = f^0(x-tv,v)$, it is clear that \eqref{e:endpointtrans} implies
\begin{equation} \label{e:endpointrho} 
\|\rho(f^0)\|_{L^{\frac{d+1}{d}}_tL^{\frac{d+1}{d-1}}_x} \lesssim \|f^0\|_{L^{\frac{d+1}{d}}_{x,v}},
\end{equation}
where $\rho(f^0)$ is the macroscopic density defined by the linear mapping 
\begin{equation*}
\rho(f^0)(t,x) = \int_{\mathbb{R}^d} f^0(x-tv,v)\, \mathrm{d}v.
\end{equation*}
Hence, by duality, \eqref{e:endpointrho} implies
\begin{equation} \label{e:endpoint}
\| \rho^* g \|_{L^{d+1}_{x,v}} \lesssim \| g \|_{L^{d+1}_t L^{\frac{d+1}{2}}_x},
\end{equation}
where the adjoint $\rho^*$ is given by
\begin{equation*}
\rho^*g(x,v) = \int_{\mathbb{R}} g(t,x+tv)\,\mathrm{d}t.
\end{equation*}
From here, the argument strongly uses ideas from the paper of Frank \emph{et al.} \cite{Franketal} concerning refined Strichartz estimates for the free Schr\"odinger propagator associated with orthonormal initial data.

Suppose $g \in \mathcal{S}(\mathbb{R} \times \mathbb{R}^d) \setminus \{0\}$ is nonnegative and such that $\widehat{g} \in C^\infty_c(\mathbb{R} \times \mathbb{R}^d)$ is also nonnegative. Here, we use $\widehat{g}$ to denote the space-time Fourier transform of $g$ given by
$$
\widehat{g}(\tau,\xi) = \int_{\mathbb{R} \times \mathbb{R}^d} g(t,x) e^{-i(t\tau + x \cdot \xi)} \, \mathrm{d}t\mathrm{d}x.
$$
In this proof, we shall also use $c$ to denote a constant depending on at most $d$, which may change from line to line.

Proceeding formally, using Fourier inversion we get
\begin{eqnarray*}
\| \rho^* g\|^{d+1}_{L^{d+1}_{x,v}} & = & \int \prod_{j=1}^{d+1} g(t_j,x+t_jv) \,\mathrm{d}\vec{t} \mathrm{d}x \mathrm{d}v \\
& = & c \int \prod_{j=1}^{d+1} \widehat{g}(\tau_j,\xi_j) \prod_{k=1}^{d+1} e^{i t_k(\tau_k + v \cdot \xi_k)}  e^{i x \cdot \sum_{\ell=1}^{d+1} \xi_\ell} \,  \mathrm{d}\vec{\tau} \mathrm{d}\vec{\xi} \mathrm{d}\vec{t} \mathrm{d}x \mathrm{d}v \\
& = & c \int \prod_{j=1}^{d+1} \widehat{g}(\tau_j,\xi_j)\, \prod_{k=1}^{d+1} \delta(\tau_k + v \cdot \xi_k) \, \,\delta\bigg(\sum_{\ell = 1}^{d+1} \xi_\ell \bigg) \mathrm{d}\vec{\tau} \mathrm{d}\vec{\xi} \mathrm{d}v \\
& = & c \int \prod_{j=1}^{d+1} \widehat{g}(-v \cdot \xi_j, \xi_j) \,\delta\bigg(\sum_{k=1}^{d+1} \xi_k\bigg) \mathrm{d}\vec{\xi}  \mathrm{d}v
\end{eqnarray*}
and hence
\begin{equation} \label{e:maincounter}
\| \rho^* g\|^{d+1}_{L^{d+1}_{x,v}} = c\int \prod_{j=1}^{d} \widehat{g}(-v \cdot \xi_j, \xi_j) \,\, \widehat{g}\bigg(v \cdot \sum_{k=1}^d \xi_k, -\sum_{\ell=1}^d \xi_\ell\bigg) \, \mathrm{d}\vec{\xi} \mathrm{d}v.
\end{equation}
We remark that by appropriately truncating the integrals in the above identities and limiting arguments, (2.4) makes sense in $[0,\infty]$ for the class of $g$ under consideration.

Define $K$ to be the $d$ by $d$ matrix whose consecutive rows are $-\xi_1, \cdots, -\xi_d$. Using the change of variables $w = Kv$, so that $w_j = -\xi_j \cdot v$ for each $1 \leq j \leq d$, we obtain
\begin{equation*} 
\| \rho^* g\|^{d+1}_{L^{d+1}_{x,v}} = c\int \prod_{j=1}^{d} \widehat{g}(w_j, \xi_j) \,\, \widehat{g}\bigg(-\sum_{k=1}^d w_k, -\sum_{\ell=1}^d \xi_\ell \bigg) \, \frac{1}{|\det K|} \,\mathrm{d}\vec{w} \mathrm{d}\vec{\xi}\,.
\end{equation*}
Writing each $\xi_j = r_j \theta_j$ in polar coordinates, we have
\begin{equation*}
|\det K| = \bigg(\prod_{j=1}^d r_j\bigg) \, |\det(\theta_1 \cdots \theta_d)|.
\end{equation*}
Since $\widehat{g}(0,0) > 0$ and $\widehat{g}$ is continuous, it follows that
\begin{equation} \label{e:final} 
\| \rho^* g\|^{d+1}_{L^{d+1}_{x,v}} \gtrsim \int_{|r| \lesssim 1} \int_{{(\mathbb{S}^{d-1})^d} } \bigg(\prod_{j=1}^{d} r_j^{d-2}\bigg) \, \frac{1}{|\det (\theta_1 \cdots \theta_d)|} \,\mathrm{d}\vec{r} \mathrm{d}\vec{\theta}. \end{equation}
For $d=1$ the radial integral is infinite, and for $d \geq 2$, 
\begin{equation*}
\int_{(\mathbb{S}^{d-1})^d} \frac{1}{|\det (\theta_1 \cdots \theta_d)|}  \,\mathrm{d}\vec{\theta} = \infty,
\end{equation*}
so the angular integral is infinite. Hence, for all $d \geq 1$ we have shown that \eqref{e:endpoint}, and consequently \eqref{e:endpointL2}, cannot hold.
 
\begin{remark}
The above argument shows that the endpoint estimate \eqref{e:endpoint} fails rather generically. For example, the space-velocity norm $\|\rho^*g\|_{d+1}$ is infinite whenever $g \in \mathcal{S}(\mathbb{R} \times \mathbb{R}^d) \setminus \{0\}$ is nonnegative and such that $\widehat{g} \in C^\infty_c(\mathbb{R} \times \mathbb{R}^d)$ is also nonnegative. 
\end{remark}

\section{A multilinear approach to the non-endpoint cases}
Fix $\sigma > 1$. In this section, the notation $\lesssim$ allows, in addition, the implicit constant to depend on $\sigma$. 

We shall prove 
\begin{equation} \label{e:p>1}
\| \rho^*g \|_{L^{\sigma(d+1)}_{x,v}} \lesssim \|g\|_{L^{q(\sigma)}_tL^{\frac{(d+1)\sigma}{2}}_x},
\end{equation}
for all $g \in L^{q(\sigma)}_tL^{\frac{(d+1)\sigma}{2}}_x$, where the exponent $q(\sigma)$ satisfies
\begin{equation*}
\frac{1}{q(\sigma)} + \frac{d}{(d+1)\sigma} = 1.
\end{equation*}

Using the invariance under transformations in \eqref{e:transformations}, to prove
the full range of non-endpoint Strichartz estimates, it suffices
to consider $(q,p,a)$ satisfying 
\begin{equation} \label{e:r1scaling}
q > a, \qquad p \geq a, \qquad \frac{2}{q} = d\bigg(1-\frac{1}{p}\bigg), \qquad  \frac{1}{a} = \frac{1}{2}\bigg(1 + \frac{1}{p}\bigg)
\end{equation}
and show that \eqref{e:Strichartz} holds with $r=1$, or equivalently, that
\begin{equation*}
\| \rho(f^0) \|_{L^q_tL^p_x} \lesssim \|f^0\|_{L^a_{x,v}}
\end{equation*}
holds for all $f^0 \in L^a_{x,v}$. By duality, this is equivalent to
\begin{equation} \label{e:almostthere}
\| \rho^*g \|_{L^{a'}_{x,v}} \lesssim \|g\|_{L^{q'}_tL^{p'}_x}
\end{equation}
for all $g \in L^{q'}_tL^{p'}_x$. Note that \eqref{e:r1scaling} implies that $a' = 2p'$ and $\frac{1}{q'} + \frac{d}{a'}  =  1$, in which case \eqref{e:almostthere} reads
\begin{equation*} 
\| \rho^*g \|_{L^{a'}_{x,v}} \lesssim \|g\|_{L^{q'}_tL^{\frac{a'}{2}}_x},
\end{equation*}
and the condition $q > a$ is equivalent to $a' > d+1$. Therefore, \eqref{e:p>1} with $\sigma = \frac{a'}{d+1} > 1$ implies the full range of non-endpoint Strichartz estimates \eqref{e:Strichartz}.

\begin{proof}[Proof of \eqref{e:p>1}]
Without loss of generality, suppose $g$ is nonnegative. By multiplying out and using Minkowski's integral inequality, we get
\begin{equation*}
\| \rho^*g \|_{L^{\sigma(d+1)}_{x,v}}^{d+1} = \bigg(\int \bigg( \int \prod_{j=1}^{d+1} g(t_j,x+t_j v) \, \mathrm{d}\vec{t} \bigg)^\sigma \, \mathrm{d}x \mathrm{d}v \bigg)^{1/\sigma} \leq  \int \bigg( \int \prod_{j=1}^{d+1} g(t_j,x+t_j v)^\sigma \, \mathrm{d}x \mathrm{d}v \bigg)^{1/\sigma} \, \mathrm{d}\vec{t}.
\end{equation*}

Now fix $t_1, \ldots, t_{d+1}$ and consider the multilinear form
$$
\int \prod_{j=1}^{d+1} g_j(t_j,x+t_j v) \, \mathrm{d}x \mathrm{d}v.
$$
A straightforward estimate via the change of variables $(x,v) \mapsto (x+t_iv,x+t_jv)$ gives
$$
\int \prod_{j=1}^{d+1} g_j(t_j,x+t_j v) \, \mathrm{d}x \mathrm{d}v \lesssim \frac{1}{|t_i - t_j|^d} \|g_i(t_i,\cdot)\|_{L^1_x} \|g_j(t_j,\cdot)\|_{L^1_x} \prod_{k \neq i,j} \|g_k(t_k, \cdot)\|_{L^\infty_x}
$$
for each $1 \leq i < j \leq d$. A multilinear interpolation argument yields
$$
\int \prod_{j=1}^{d+1} g_j(t_j,x+t_j v) \, \mathrm{d}x \mathrm{d}v \lesssim \prod_{1 \leq i < j \leq d} |t_i - t_j|^{-\frac{2}{d+1}} \,\, \prod_{k=1}^{d+1} \|g_k(t_k,\cdot)\|_{L^{\frac{d+1}{2}}_x}. 
$$
Applying this with $g^\sigma$ for each $g_j$ we get
\begin{eqnarray*}
\| \rho^*g \|_{L^{\sigma(d+1)}_{x,v}}^{d+1} & \lesssim & \int_{\mathbb{R}^{d+1}} \prod_{1 \leq i < j \leq d} |t_i - t_j|^{-\frac{2}{(d+1)\sigma}} \,\, \prod_{k=1}^{d+1} \|g(t_k,\cdot)\|_{L^{\frac{(d+1)\sigma}{2}}_x} \, \mathrm{d}\vec{t} \\
& \lesssim & \|g\|_{L^{q(\sigma)}_tL^{\frac{(d+1)\sigma}{2}}_x}^{d+1}\,.
\end{eqnarray*}
The last inequality is a consequence of the multilinear Hardy--Littlewood--Sobolev inequality due to Christ \cite{Christ}. 
\end{proof}

\begin{remark} Certain replacements for the endpoint are already known. For example, with $q=a=2$, in \cite{KeelTao}, Keel and Tao obtain a substitute for \eqref{e:endpointL2} with the
$L^{p}_xL^{r}_v$ norm (where $p = \frac{2d}{d-1}$ and $r = \frac{2d}{d+1}$) replaced by that of a certain real
interpolation space which is between $L^{p,1}_xL^{r,1}_v$ and
$L^{p,\infty}_xL^{r,\infty}_v$. See also work of Ovcharov \cite{Ovcharov} where a different substitute bound was given for velocities $v$ belonging to a bounded subset of $\mathbb{R}^d$. 
\end{remark}

\begin{acknowledgement}
We would like to thank Jos\'e Ca\~nizo for his contributions during the early stages of this project. This work was supported by the European Research Council [grant
number 307617] (Bennett); the Engineering and Physical Sciences Research Council [grant numbers EP/J021490/1, EP/J01155X/1] (Bez, Guti\'errez); National Research Foundation of Korea [grant number 2012008373] (Lee).
\end{acknowledgement}


\begin{thebibliography}{99}

\bibitem{CastellaPerthame} F. Castella, B. Perthame, \textit{Estimations de Strichartz pour les \`equations de transport cin\'etique}, C. R.
Acad. Sci. Paris S\'er. I Math. \textbf{332} (1996), 535--540.

\bibitem{Christ} M. Christ, \textit{On the restriction of the Fourier transform to curves: endpoint results and the degenerate case}, Trans. Amer. Math. Soc, \textbf{287} (1985), 223--238.

\bibitem{Franketal} R. Frank, M. Lewin, E. H. Lieb, R. Seiringer, \textit{A Strichartz inequality for orthonormal functions}, arXiv:1306.1309.

\bibitem{GuoPeng} Z. Guo, L. Peng, \textit{Endpoint Strichartz estimate for the kinetic transport
equation in one dimension}, C. R. Math. Acad. Sci. Paris, \textbf{345} (2007), 253--256.

\bibitem{KeelTao} M. Keel, T. Tao, \textit{Endpoint Strichartz estimates}, Amer. J. Math., \textbf{120} (1998), 955--980.

\bibitem{MontSmith} S. J. Montgomery-Smith, \textit{Time decay for the bounded mean oscillation of solutions of the Schr\"odinger and wave equations}, Duke Math. J., \textbf{91} (1998), 393--408.

\bibitem{OvcharovKakeya} E. Ovcharov, \textit{Counterexamples to Strichartz estimates for the kinetic
transport equation based on Besicovitch sets}, Nonlinear Analysis: Theory, Methods \& Applications
\textbf{74} (2011), 2515--2522.

\bibitem{Ovcharov} E. Ovcharov, \textit{Strichartz estimates for the kinetic transport equation},
SIAM J. Math. Anal., \textbf{43} (2011), 1282--1310.
\end{thebibliography}
\end{document}